\documentclass[12pt]{amsart}

\def\RPone{Require}
\def\RPtwo{Pack}
\def\RPthree{age}
\expandafter\let\expandafter\LoadStyle\csname\RPone\RPtwo\RPthree\endcsname
\LoadStyle{cite}
\LoadStyle[a4paper,left=2cm,right=2cm,top=2cm,bottom=2cm]{geometry}
\LoadStyle{amsmath,amssymb,amsthm,mathtools}
\mathtoolsset{showonlyrefs=true}
\LoadStyle[T1]{fontenc}
\LoadStyle{lmodern}
\LoadStyle[colorlinks=true,linkcolor=blue,citecolor=blue,urlcolor=blue]{hyperref}
\theoremstyle{plain}
\newtheorem{theorem}{Theorem}[section]

\newtheorem{proposition}[theorem]{Proposition}
\newtheorem{corollary}[theorem]{Corollary}
\theoremstyle{definition}
\newtheorem{definition}[theorem]{Definition}
\newtheorem{example}[theorem]{Example}
\newtheorem{remark}[theorem]{Remark}
\newtheorem{notation}[theorem]{Notation}

\newcommand{\bbN}{\mathbb N}

\newcommand{\calI}{\mathcal I}
\newcommand{\calP}{\mathcal P}
\newcommand{\Id}{\operatorname{Id}}
\newcommand{\Spec}{\operatorname{Spec}}
\newcommand{\SpecPhi}{\operatorname{Spec}_{\!\Phi}}
\newcommand{\Phiop}{\Phi^{\!*}}

\newcommand{\m}{\mathfrak m}
\newcommand{\Zar}{\mathrm{Zar}}
\newcommand{\Ker}{\operatorname{Ker}}
\newcommand{\Ann}{\operatorname{Ann}}

\newcommand{\Sob}{\operatorname{Sob}}
\newcommand{\Vaff}{V_{\mathrm{aff}}}
\newcommand{\VSpec}{V_{\operatorname{Spec}}}
\newcommand{\setprod}{\mathbin{\odot}}
\def\DD{D\kern-.7em\raise0.4ex\hbox{\char '55}\kern.33em}

\title[Observation nuclei and infinitesimal recovery]{Observation Nuclei in Affine Geometry:\\ Symbolic Powers, Fat Points, and Infinitesimal Recovery}

\author{\DD\d{\u a}ng V\~o Ph\'uc}
\address{Department of Mathematics, FPT University, An Phu Thinh New Urban Area, Quy Nhon, Vietnam}
\email{dangphuc150488@gmail.com}
\thanks{ORCID: \url{https://orcid.org/0000-0002-6885-3996}}


\begin{document}

\begin{abstract}
This paper integrates point-dependent infinitesimal data into affine geometry through the established frameworks of quantale nuclei and finitary weak ideal systems. For a subset \(X\subseteq k^n\) and a bounded function \(\nu:X\to\mathbb N_{>0}\), we study the observation nucleus
\[
\Phiop_{X,\nu}(A)=\bigcap_{a\in X}\bigl(A+\m_a^{\nu(a)}\bigr).
\]
The construction organizes several classical interfaces within one point-dependent operator: closed observations recover symbolic powers through the Zariski--Nagata theorem, and finite observations restrict on ideals to translation by a fat-point ideal through the Chinese remainder theorem. The principal uniform identity is the exact relative Nullstellensatz
\[
\sqrt{\Phiop_{X,\nu}(I)}=\calI\bigl(X\cap\Vaff(I)\bigr),
\]
which separates the reduced geometry determined by the observation set from the retained infinitesimal multiplicities determined by \(\nu\). This identity yields a density criterion for fixed prime ideals and identifies their space with the sobrification of the observed point set. For constant full-affine observation order, the tower has common radical \(\sqrt I\), and the nilpotent Nullstellensatz of Eisenbud and Hochster implies finite-stage stabilization: for every ideal \(I\), there exists \(e=e(I)\geq1\) such that \(\Phiop_q(I)=I\) for all \(q\geq e\). The observation formalism therefore gives a simultaneous closed-point congruence interpretation of this classical theorem. Finally, the full affine nucleus is not a fixed-ideal translation in any positive dimension, and an explicit unbounded example shows that boundedness cannot be omitted from the general radical formula.
\end{abstract}

\maketitle

\section{Introduction}\label{sec:introduction}

Approximation constructions in algebra have been formulated through upper approximations attached to proximity, descriptive data, and tolerance relations. Peters~\cite{Peters2013} provided a general introduction to near sets and descriptive proximity, emphasizing the role of nearness relations and feature-based comparison in proximity spaces. Approximately rings and approximately ideals in proximal relator spaces were introduced by \.{I}nan \cite{Inan2019}, and approximately prime rings and ideals were subsequently studied by Almahariq, Peters, and Vergili \cite{APV2025}. Once an upper approximation is assumed to be extensive, monotone, idempotent, and compatible with the ring operations, its ideal-theoretic consequences enter established order-theoretic and multiplicative frameworks.

The first framework is the theory of quantales. Rosenthal's monograph treats nuclei on quantales, and Elliott places quantic nuclei, star and semistar operations, semiprime operations, ideal systems, and module systems in a common closure-theoretic setting \cite{Rosenthal1990,Elliott2015,Elliott2019}. The second framework is multiplicative ideal theory. On the ideal lattice of a commutative ring, a finite-type nucleus determines a finitary weak ideal system in the sense developed by Halter-Koch \cite{HalterKoch1998,HalterKoch2001,HalterKoch2011,HalterKoch2025}. Prime separation and radical representation are inherited from this theory. Spectrality of prime spectra attached to finitary weak ideal systems is likewise known; see Juett \cite[Theorem~5.9]{Juett2012} and Finocchiaro--Yazdonov \cite{FinocchiaroYazdonov2026}. We retain a direct proof in Section~\ref{sec:spectral} only to establish the notation and the precise topological form used in the observation construction.

While the abstract algebraic properties of these operators are governed by ideal systems, their concrete geometric behavior requires specific localized constructions. In this paper, we address the problem of integrating point-dependent infinitesimal data into a uniform affine geometric framework. Our approach utilizes the observation nucleus $\Phiop_{X,\nu}$, defined via simultaneous congruence saturations over a subset of points $X \subseteq k^n$ equipped with a bounded order function $\nu$. The geometric input comes from infinitesimal neighborhoods of closed points. For a point $a\in k^n$, the quotient $k[x_1,\ldots,x_n]/\m_a^q$ is the coordinate algebra of the $(q-1)$-st infinitesimal neighborhood of $a$. Its classes record shifted Taylor or Hasse data through total degree less than $q$. This local algebra must be distinguished from the Hasse--Schmidt algebras representing schemes of truncated arcs; see Vojta \cite{Vojta2007} and Musta\c{t}\u{a} \cite{Mustata2001}. It is also related, but not identical, to the local jet closures introduced by de Fernex, Ein, and Ishii \cite{deFernexEinIshii2018} and developed further by Chen and Zuo \cite{ChenZuo2025}. Their closures are defined through equality of local jet schemes or jet supports, whereas the factors used here are elementary congruence saturations $I+\m_a^q$. Known coincidences in the homogeneous case are specified in Remark~\ref{rem:hasse-jet}.

This distinction gives a concrete computational advantage for congruence membership. Fix a point $a$ and an order $q$, and write $I=(f_1,\ldots,f_s)$. In the shifted coordinates $z_i=x_i-a_i$, the algebra
\[
R/\m_a^q\cong k[z_1,\ldots,z_n]/(z_1,\ldots,z_n)^q
\]
has basis the monomials $z^\alpha$ with $|\alpha|<q$, and hence dimension $\binom{n+q-1}{n}$. The image of $I$ is the $k$-linear span of the truncated shifted expansions of $z^\beta f_i(a+z)$ with $|\beta|<q$. Indeed, the coefficient multiplying each generator can itself be reduced to a linear combination of those monomials. Thus $g\in I+\m_a^q$ is decided by linear algebra on its shifted Taylor coefficient vector, using Hasse coefficients in arbitrary characteristic; this is a finite quotient-algebra computation of the kind discussed in \cite{CoxLittleOShea2015}. The congruence saturation therefore bypasses the generally nonlinear equations in unknown arc coefficients that arise in direct jet-scheme constructions \cite{Vojta2007}. The advantage concerns these local membership tests; an arbitrary infinite observation set still requires a separate argument to determine the global intersection.

The proposed methodology unifies several classical constructions into a single point-dependent operator on arbitrary subsets. Symbolic powers and fat-point schemes provide two primary specializations. If $k$ is algebraically closed, $P\subseteq k[x_1,\ldots,x_n]$ is prime, and the observation set is the affine variety $\Vaff(P)$ with constant order $q$, the Zariski--Nagata theorem gives $\Phiop_{\Vaff(P),q}(0) = \bigcap_{a\in\Vaff(P)}\m_a^q = P^{(q)}$; see \cite{DaoEtAl2018}. A contemporary differential-operator characterization of symbolic powers over arbitrary fields is given by De Stefani, Grifo, and Jeffries \cite{DeStefaniGrifoJeffries2025}. If the observation set is finite, pairwise comaximality and the Chinese remainder theorem \cite{AtiyahMacdonald1969} show that the restriction to ideals is translation by the corresponding fat-point ideal. These identifications delimit the novelty of the framework: the present construction does not replace symbolic-power, fat-point, or jet-closure theory, but rather assembles their local congruence patterns to study their global affine properties.

The main results of this paper are organized around the exactness, topological classification, and structural recovery of this observation nucleus. Specifically:

First, we establish the exact relative Nullstellensatz for a bounded observation function $\nu:X\to\mathbb N_{>0}$. We prove that $\sqrt{\Phiop_{X,\nu}(I)} = \calI\bigl(X\cap\Vaff(I)\bigr)$. This result guarantees that the observation set controls the reduced geometry, while the values of $\nu$ precisely control the retained infinitesimal multiplicity.

Second, based on this exact Nullstellensatz, we characterize the fixed-prime spectrum. We show that a prime ideal $P$ is fixed precisely when $X\cap\Vaff(P)$ is Zariski dense in $\Vaff(P)$. The corresponding fixed-prime space is canonically identified with the sobrification of $X$ equipped with its induced Zariski topology. The underlying topological construction is the established sobrification of a subspace of a sober space \cite[Theorem~4.20]{BouacidaEtAl2003}; see also \cite[Sections~2.2--2.3]{Schwartz2023}. Here it is realized by the fixed-prime locus of the observation nucleus.

Third, for the full affine observation set and constant order $q$, the operators form a descending tower
\[
\sqrt I=\Phiop_1(I)\supseteq\Phiop_2(I)\supseteq\cdots\supseteq I
\]
with common radical $\sqrt I$. The quotient $\Phiop_q(I)/I$ is the intersection of the $q$-th powers of all closed-point maximal ideals of $R/I$. Eisenbud and Hochster's Nullstellensatz with nilpotents \cite{EisenbudHochster1979} then supplies a uniform exponent $e=e(I)$ for which $\Phiop_q(I)=I$ for every $q\geq e$. Thus finite-stage recovery is classical; the role of the present construction is to place it inside the observation-nucleus framework and to compare it with the relative, point-dependent, and nontranslation phenomena.

Finally, we demonstrate the non-triviality and sharpness of the model. We prove that the full affine nucleus is not a fixed-ideal translation in any positive dimension, and we provide an explicit unbounded example to prove that the uniform boundedness hypothesis cannot be omitted from the general radical formula.

The paper is organized as follows. Section~\ref{sec:synthesis} records the quantale and weak-ideal-system synthesis. Section~\ref{sec:spectral} states the fixed-prime spectral theory in the notation used later. Section~\ref{sec:observation} constructs observation nuclei and establishes their relations with infinitesimal neighborhoods, symbolic powers, and local jet closures. Section~\ref{sec:nullstellensatz} proves the exact relative Nullstellensatz, the density criterion, and the sobrification description. Finally, Section~\ref{sec:tower} develops the closed-point power interpretation, the Eisenbud--Hochster finite-stage recovery theorem, the all-dimensional nontranslation theorem, the unbounded-order counterexample, and the finite fat-point formula.

\section{Nuclei and Finitary Weak Ideal Systems}\label{sec:synthesis}

This section isolates the algebraic substrate used throughout the paper. The goal is to pass from a closure on arbitrary subsets to a nucleus on the ideal quantale and then invoke the standard prime-separation theorem for finitary weak ideal systems.

\begin{notation}\label{not:subset-operations}
Let \(R\) be a commutative ring with unity. For subsets \(A,B\subseteq R\), define
\[
A+B=\{a+b:a\in A,\ b\in B\},
\qquad
A\setprod B=\{ab:a\in A,\ b\in B\}.
\]
For arbitrary subsets, \(A\setprod B\) denotes pointwise multiplication. For ideals, juxtaposition \(IJ\) denotes ordinary ideal multiplication, whose elements are finite sums of products. Let \(\Id(R)\) be the complete lattice of ideals of \(R\), ordered by inclusion and equipped with ordinary ideal multiplication.
\end{notation}

\begin{definition}[Algebra-compatible subset closure]\label{def:algebra-compatible}
Let \(R\) be a commutative ring with unity. An operator \(\Phiop:\calP(R)\to\calP(R)\) is an algebra-compatible subset closure if, for all subsets \(A,B\subseteq R\) and all \(r\in R\),
\[
A\subseteq\Phiop(A),
\qquad
A\subseteq B\Longrightarrow\Phiop(A)\subseteq\Phiop(B),
\qquad
\Phiop(\Phiop(A))=\Phiop(A),
\tag{C1--C3}
\]
and
\[
\Phiop(A)+\Phiop(B)\subseteq\Phiop(A+B),
\qquad
r\Phiop(A)\subseteq\Phiop(rA).
\tag{C4}
\]
\end{definition}

\begin{definition}[Quantale nucleus]\label{def:quantale-nucleus}
Let \(Q\) be a unital quantale with multiplication \(\cdot\). An operator \(j:Q\to Q\) is a nucleus if it is extensive, monotone, idempotent, and satisfies
\[
j(x)\cdot j(y)\leq j(x\cdot y)
\]
for all \(x,y\in Q\).
\end{definition}

\begin{proposition}\label{prop:quantale-synthesis}
Let \(R\) be a commutative ring with unity, and let \(\Phiop\) be an algebra-compatible subset closure. Then \(\Phiop\) is a nucleus on the additive powerset quantale \((\calP(R),\cup,+)\) and on the multiplicative powerset quantale \((\calP(R),\cup,\setprod)\). Moreover, the assignment
\[
\rho:\Id(R)\longrightarrow\Id(R),
\qquad
\rho(I)=\Phiop(I),
\]
is a nucleus on the ideal quantale \(\Id(R)\).
\end{proposition}

\begin{proof}
The additive nucleus inequality is the first inclusion in \textup{(C4)}. We first prove the multiplicative inequality for pointwise products. Let \(x\in\Phiop(A)\) and \(y\in\Phiop(B)\). For every \(b\in B\), scalar compatibility gives
\[
xb\in b\Phiop(A)\subseteq\Phiop(bA)\subseteq\Phiop(A\setprod B).
\]
Hence
\[
xB\subseteq\Phiop(A\setprod B).
\]
Applying scalar compatibility to \(y\in\Phiop(B)\), with scalar \(x\), gives
\[
xy\in x\Phiop(B)\subseteq\Phiop(xB).
\]
Monotonicity and idempotence now yield
\[
\Phiop(xB)
\subseteq
\Phiop\bigl(\Phiop(A\setprod B)\bigr)
=
\Phiop(A\setprod B).
\]
Thus
\[
\Phiop(A)\setprod\Phiop(B)
\subseteq
\Phiop(A\setprod B),
\]
which is the multiplicative powerset-nucleus inequality.

Let \(I\lhd R\). We verify that \(\Phiop(I)\) is an ideal. Since \((-1)I=I\), scalar compatibility gives
\[
-\Phiop(I)\subseteq\Phiop(-I)=\Phiop(I).
\]
If \(u,v\in\Phiop(I)\), then
\[
u-v\in\Phiop(I)+\Phiop(I)
\subseteq\Phiop(I+I)=\Phiop(I).
\]
Thus \(\Phiop(I)\) is an additive subgroup. If \(r\in R\) and \(u\in\Phiop(I)\), then
\[
ru\in r\Phiop(I)
\subseteq\Phiop(rI)
\subseteq\Phiop(I),
\]
because \(rI\subseteq I\). Hence \(\Phiop(I)\) is an ideal.

The restrictions of extensivity, monotonicity, and idempotence make \(\rho\) a closure operator on \(\Id(R)\). Let \(I,J\lhd R\). The pointwise-product inequality proves that every product \(xy\), with \(x\in\Phiop(I)\) and \(y\in\Phiop(J)\), belongs to
\[
\Phiop(I\setprod J)\subseteq\Phiop(IJ),
\]
because \(I\setprod J\subseteq IJ\). Every element of the ordinary ideal product \(\rho(I)\rho(J)\) is a finite sum
\[
\sum_{i=1}^{t}x_i y_i
\qquad
\bigl(x_i\in\Phiop(I),\ y_i\in\Phiop(J)\bigr).
\]
Each summand belongs to the ideal \(\Phiop(IJ)\), and therefore the finite sum belongs to \(\Phiop(IJ)\). Hence
\[
\rho(I)\rho(J)\subseteq\rho(IJ).
\]
Therefore \(\rho\) is a nucleus on the ideal quantale.
\end{proof}

\begin{remark}\label{rem:quantale-context}
Proposition~\ref{prop:quantale-synthesis} is an instance of the standard nucleus construction for quantales; see \cite{Rosenthal1990,Elliott2015,Elliott2019}. The distinction between pointwise subset multiplication and ordinary ideal multiplication is essential in the final finite-sum step.
\end{remark}

\begin{definition}[Finite type]\label{def:finite-type}
Let \(R\) be a commutative ring with unity, and let \(\rho\) be a nucleus on \(\Id(R)\). The nucleus \(\rho\) is of finite type if
\[
\rho(I)=\bigcup_{\substack{J\subseteq I\\J\text{ finitely generated}}}\rho(J)
\]
for every ideal \(I\lhd R\).
\end{definition}

\begin{remark}\label{rem:weak-ideal-system}
Let
\[
\rho:\Id(R)\longrightarrow\Id(R)
\]
be a finite-type nucleus. For every subset \(X\subseteq R\), let
\((X)\) denote the ordinary ideal generated by \(X\), and define
\[
X^r:=\rho((X)).
\]
Then \(r\) is the finitary weak ideal system induced by \(\rho\), and
its \(r\)-ideals are precisely the \(\rho\)-fixed ordinary ideals.
This is the precise sense in which a finite-type nucleus on
\(\Id(R)\) gives the ideal-lattice realization of a finitary weak
ideal system.
\end{remark}

\begin{definition}[Fixed ideals and fixed primes]\label{def:fixed-prime}
Let \(R\) be a commutative ring with unity, and let \(\rho\) be a nucleus on \(\Id(R)\). An ideal \(I\lhd R\) is a \(\Phi\)-ideal if \(\rho(I)=I\). A prime ideal \(P\in\Spec(R)\) is a \(\Phi\)-prime if \(\rho(P)=P\). Define
\[
\SpecPhi(R)=\{P\in\Spec(R):\rho(P)=P\}.
\]
\end{definition}

\begin{proposition}\label{prop:collapse-prime}
Let \(R\) be a commutative ring with unity, and let \(\Phiop\) be an algebra-compatible subset closure. Let \(P\subsetneq R\) be an additive subgroup satisfying \(RP\subseteq\Phiop(P)\). If
\[
xy\in\Phiop(P)\Longrightarrow x\in P\text{ or }y\in P
\]
for all \(x,y\in R\), then \(P\) is a \(\Phi\)-prime ideal.
\end{proposition}

\begin{proof}
Assume that \(1\in P\). Then \(R=RP\subseteq\Phiop(P)\). Choose \(a\in R\setminus P\). Since \(a^2\in\Phiop(P)=R\), the displayed primality condition gives \(a\in P\), a contradiction. Hence \(1\notin P\).

Let \(x\in\Phiop(P)\). Since \(x=x\cdot1\), the primality condition gives \(x\in P\) or \(1\in P\). The second alternative is impossible, so \(x\in P\). Therefore \(\Phiop(P)\subseteq P\), while extensivity gives \(P\subseteq\Phiop(P)\). Thus \(\Phiop(P)=P\). Since \(P\) is an additive subgroup and \(RP\subseteq\Phiop(P)=P\), the set \(P\) is an ideal. If \(xy\in P=\Phiop(P)\), the displayed condition gives \(x\in P\) or \(y\in P\). Therefore \(P\) is an ordinary prime ideal fixed by \(\rho\).
\end{proof}

\begin{theorem}\label{thm:prime-separation}
Let \(R\) be a commutative ring with unity, and let \(\rho\) be a finite-type nucleus on \(\Id(R)\). Then, for every ideal \(I\lhd R\),
\[
\sqrt{\rho(I)}
=
\bigcap_{\substack{P\in\SpecPhi(R)\\I\subseteq P}}P,
\]
where the intersection of an empty family is \(R\).
\end{theorem}

\begin{proof}
The corresponding maximal-disjoint prime theorem for finitary weak ideal
systems on cancellative commutative monoids is
\cite[Proposition~6.3, p.~214]{HalterKoch2011}. Since the present theorem
allows commutative rings with zero divisors, we give the ideal-lattice
argument.

We first prove
\[
\sqrt{\rho(I)}
\subseteq
\bigcap_{\substack{P\in\SpecPhi(R)\\ I\subseteq P}}P.
\]
Let \(x\in\sqrt{\rho(I)}\), and let \(P\in\SpecPhi(R)\) satisfy
\(I\subseteq P\). There exists a positive integer \(n\) such that
\[
x^n\in\rho(I).
\]
Since \(I\subseteq P\), monotonicity of \(\rho\) and the fixedness of \(P\)
give
\[
\rho(I)\subseteq\rho(P)=P.
\]
Thus \(x^n\in P\). Since \(P\) is an ordinary prime ideal, it follows that
\(x\in P\). Therefore \(x\) belongs to every \(\Phi\)-prime containing
\(I\), proving the first inclusion.

For the reverse inclusion, let
\[
x\notin\sqrt{\rho(I)}.
\]
Define the multiplicatively closed set
\[
T_x=\{1,x,x^2,x^3,\ldots\}.
\]
We claim that
\[
\rho(I)\cap T_x=\varnothing.
\]
Indeed, if \(x^n\in\rho(I)\) for some positive integer \(n\), then
\(x\in\sqrt{\rho(I)}\), contrary to the choice of \(x\). Moreover,
\(1\notin\rho(I)\), because \(1\in\rho(I)\) would imply
\(\rho(I)=R\) and hence \(\sqrt{\rho(I)}=R\). Thus no element of \(T_x\)
belongs to \(\rho(I)\).

Consider the partially ordered set
\[
\mathcal C
=
\left\{
J\lhd R:
\rho(J)=J,\ 
\rho(I)\subseteq J,\ 
J\cap T_x=\varnothing
\right\},
\]
ordered by inclusion. Since \(\rho(I)\) is fixed by idempotence and is
disjoint from \(T_x\), one has \(\rho(I)\in\mathcal C\). Hence
\(\mathcal C\neq\varnothing\).

Let \(\{J_\lambda\}_{\lambda\in\Lambda}\) be a nonempty chain in \(\mathcal C\), and
set
\[
J=\bigcup_{\lambda\in\Lambda}J_\lambda.
\]
Because the family is linearly ordered by inclusion, \(J\) is an ideal.
It contains \(\rho(I)\), and it is disjoint from \(T_x\).

We verify that \(J\) is fixed by \(\rho\). Extensivity gives
\[
J\subseteq\rho(J).
\]
Let \(y\in\rho(J)\). Since \(\rho\) is of finite type, there exists a
finitely generated ideal \(L\subseteq J\) such that
\[
y\in\rho(L).
\]
Write
\[
L=(u_1,\ldots,u_s).
\]
For each \(i\), choose \(\lambda_i\in\Lambda\) such that
\(u_i\in J_{\lambda_i}\). Since the family is a chain and only finitely
many indices occur, there exists \(\lambda_0\in\Lambda\) such that
\[
u_1,\ldots,u_s\in J_{\lambda_0}.
\]
Consequently,
\[
L\subseteq J_{\lambda_0}.
\]
Monotonicity and the fixedness of \(J_{\lambda_0}\) now give
\[
y\in\rho(L)
\subseteq
\rho(J_{\lambda_0})
=
J_{\lambda_0}
\subseteq J.
\]
Thus \(\rho(J)\subseteq J\), and hence \(\rho(J)=J\). Therefore
\(J\in\mathcal C\), so every nonempty chain in \(\mathcal C\) has an upper bound. The empty chain has \(\rho(I)\) as an upper bound.

By Zorn's lemma, \(\mathcal C\) has a maximal element \(P\). We prove that
\(P\) is prime. Suppose that
\[
ab\in P
\qquad\text{and}\qquad
a,b\notin P.
\]
Since \(a\notin P\), the fixed ideal
\[
\rho(P+(a))
\]
strictly contains \(P\). If it were disjoint from \(T_x\), it would belong
to \(\mathcal C\), contradicting the maximality of \(P\). Hence there
exists an integer \(m\geq0\) such that
\[
x^m\in\rho(P+(a)).
\]
Similarly, there exists an integer \(n\geq0\) such that
\[
x^n\in\rho(P+(b)).
\]

The nucleus inequality gives
\[
\rho(P+(a))\,\rho(P+(b))
\subseteq
\rho\bigl((P+(a))(P+(b))\bigr).
\]
Since \(P\) is an ideal and \(ab\in P\),
\[
(P+(a))(P+(b))
=
P^2+aP+bP+(ab)
\subseteq P.
\]
Therefore, by monotonicity and the fixedness of \(P\),
\[
\rho\bigl((P+(a))(P+(b))\bigr)
\subseteq
\rho(P)
=
P.
\]
It follows that
\[
x^{m+n}=x^m x^n\in P.
\]
But \(x^{m+n}\in T_x\), contradicting \(P\cap T_x=\varnothing\).
Therefore \(P\) is prime.

By construction, \(\rho(P)=P\), so \(P\in\SpecPhi(R)\). Moreover,
\[
I\subseteq\rho(I)\subseteq P,
\]
whereas \(x\notin P\), because \(x\in T_x\) and \(P\cap T_x=\varnothing\).
Thus \(x\) does not belong to the intersection of all \(\Phi\)-prime
ideals containing \(I\). We have proved
\[
\bigcap_{\substack{P\in\SpecPhi(R)\\ I\subseteq P}}P
\subseteq
\sqrt{\rho(I)}.
\]
Combining the two inclusions gives the asserted equality.
\end{proof}

\section{Spectral Geometry of Fixed Prime Ideals}\label{sec:spectral}

The algebraic synthesis permits the geometry to be formulated on the fixed-prime locus. Spectrality of prime spectra associated with finitary weak ideal systems on commutative monoids with zero was established by Juett \cite[Theorem~5.9]{Juett2012}, and a monoid-theoretic finite-type formulation appears in Finocchiaro--Yazdonov \cite{FinocchiaroYazdonov2026}. We give a direct proof in the present notation because its quasi-compact principal opens and generic-point construction are used repeatedly below. We also record how comparison of nuclei reverses the inclusion of spectra, which separates the effects of adding observers from those of increasing observation order.

\begin{definition}[Fixed-prime topology]\label{def:fixed-prime-topology}
Let \(R\) be a commutative ring with unity, and let \(\rho\) be a nucleus on \(\Id(R)\). For an ideal \(I\lhd R\), define
\[
V_{\Phi}(I)=\{P\in\SpecPhi(R):I\subseteq P\}.
\]
For \(f\in R\), define
\[
D_{\Phi}(f)=\SpecPhi(R)\setminus V_{\Phi}((f)).
\]
Equip \(\SpecPhi(R)\) with the topology whose closed sets are the sets \(V_{\Phi}(I)\).
\end{definition}

\begin{proposition}\label{prop:subspace-realization}
Let \(R\) be a commutative ring with unity, and let \(\rho\) be a nucleus on \(\Id(R)\). Then the topology in Definition~\ref{def:fixed-prime-topology} is the subspace topology inherited from \(\Spec(R)\). Moreover, the sets \(D_{\Phi}(f)\), with \(f\in R\), form a basis and satisfy
\[
D_{\Phi}(f)\cap D_{\Phi}(g)=D_{\Phi}(fg).
\]
\end{proposition}

\begin{proof}
For every ideal \(I\),
\[
V_{\Phi}(I)=\VSpec(I)\cap\SpecPhi(R).
\]
Hence the displayed sets are exactly the intersections of classical spectral Zariski closed sets with \(\SpecPhi(R)\), proving the subspace assertion \cite[Chapter~1]{AtiyahMacdonald1969}. The complements of the classical principal closed sets form a basis of \(\Spec(R)\), so their intersections with the subspace form a basis of \(\SpecPhi(R)\). Since every point of \(\SpecPhi(R)\) is an ordinary prime ideal,
\[
fg\notin P
\Longleftrightarrow f\notin P\text{ and }g\notin P.
\]
This equivalence gives the identity for principal opens.
\end{proof}

\begin{theorem}\label{thm:spectrality}
Let \(R\) be a commutative ring with unity, and let \(\rho\) be a finite-type nucleus on \(\Id(R)\). Then \(\SpecPhi(R)\) is a spectral space. More precisely, it is quasi-compact and \(T_0\), the principal opens \(D_{\Phi}(f)\) form a basis of quasi-compact opens stable under finite intersections, and every nonempty irreducible closed subset has a unique generic point.
\end{theorem}

\begin{proof}
We first prove that every principal open is quasi-compact. Let \(f\in R\), and suppose that
\[
D_{\Phi}(f)\subseteq\bigcup_{\lambda\in\Lambda}D_{\Phi}(g_{\lambda}).
\tag{3.1}\label{eq:cover}
\]
Let
\[
J=(g_{\lambda}:\lambda\in\Lambda).
\]
If a \(\Phi\)-prime \(P\) contains \(J\), then \(g_{\lambda}\in P\) for every \(\lambda\), so
\[
P\notin\bigcup_{\lambda\in\Lambda}D_{\Phi}(g_{\lambda}).
\]
The inclusion \eqref{eq:cover} therefore implies \(P\notin D_{\Phi}(f)\), and hence \(f\in P\). Thus \(f\) belongs to every \(\Phi\)-prime containing \(J\). Theorem~\ref{thm:prime-separation} gives
\[
f\in\sqrt{\rho(J)}.
\]
Consequently, there is a positive integer \(N\) such that
\[
f^N\in\rho(J).
\]
Because \(\rho\) is of finite type, there exists a finitely generated ideal \(K\subseteq J\) such that
\[
f^N\in\rho(K).
\]
Choose generators \(h_1,\ldots,h_t\) of \(K\). Since each \(h_i\) belongs to \(J\), it is a finite \(R\)-linear combination of elements among the \(g_{\lambda}\). Only finitely many indices occur in all these expressions. Let \(\Lambda_0\subseteq\Lambda\) be their union, and set
\[
J_0=(g_{\lambda}:\lambda\in\Lambda_0).
\]
Then \(K\subseteq J_0\), so monotonicity gives
\[
f^N\in\rho(K)\subseteq\rho(J_0).
\]
Let \(P\in\SpecPhi(R)\) contain \(J_0\). Since \(P\) is fixed,
\[
\rho(J_0)\subseteq\rho(P)=P.
\]
Therefore \(f^N\in P\), and ordinary primality gives \(f\in P\). We have shown that
\[
V_{\Phi}(J_0)\subseteq V_{\Phi}((f)).
\]
Taking complements yields
\[
D_{\Phi}(f)
\subseteq
\bigcup_{\lambda\in\Lambda_0}D_{\Phi}(g_{\lambda}).
\]
Thus \(D_{\Phi}(f)\) is quasi-compact. Since
\[
\SpecPhi(R)=D_{\Phi}(1),
\]
the whole space is quasi-compact. Proposition~\ref{prop:subspace-realization} gives a basis of principal opens, and the identity
\[
D_{\Phi}(f)\cap D_{\Phi}(g)=D_{\Phi}(fg)
\]
shows that this basis is stable under finite intersections.

We next prove the \(T_0\) property. Let \(P,Q\in\SpecPhi(R)\) with \(P\neq Q\). At least one of the inclusions \(P\subseteq Q\) and \(Q\subseteq P\) fails. Suppose that \(P\nsubseteq Q\), and choose \(f\in P\setminus Q\). Then
\[
Q\in D_{\Phi}(f)
\qquad\text{and}\qquad
P\notin D_{\Phi}(f).
\]
Hence the two points are topologically distinguishable, proving that the space is \(T_0\).

It remains to prove sobriety. Let \(C\) be a nonempty irreducible closed subset of \(\SpecPhi(R)\). Define
\[
P_C=\bigcap_{Q\in C}Q.
\]
The ideal \(P_C\) is proper because \(C\neq\varnothing\) and no \(Q\in C\) contains \(1\). It is fixed: extensivity gives \(P_C\subseteq\rho(P_C)\), while monotonicity gives
\[
\rho(P_C)
\subseteq
\bigcap_{Q\in C}\rho(Q)
=
\bigcap_{Q\in C}Q
=P_C.
\]
We prove that \(P_C\) is prime. Let \(a,b\in R\) satisfy \(ab\in P_C\). Then \(ab\in Q\) for every \(Q\in C\), so ordinary primality of \(Q\) gives
\[
C\subseteq V_{\Phi}((a))\cup V_{\Phi}((b)).
\]
Since \(C\) is irreducible and both sets on the right are closed, either
\[
C\subseteq V_{\Phi}((a))
\qquad\text{or}\qquad
C\subseteq V_{\Phi}((b)).
\]
In the first case, \(a\in Q\) for every \(Q\in C\), so \(a\in P_C\). In the second case, \(b\in P_C\). Hence \(P_C\) is prime and fixed, so \(P_C\in\SpecPhi(R)\).

Write \(C=V_{\Phi}(I)\) for an ideal \(I\). Every \(Q\in C\) contains \(I\), so \(I\subseteq P_C\). It follows that
\[
V_{\Phi}(P_C)\subseteq V_{\Phi}(I)=C.
\]
Conversely, every \(Q\in C\) contains the intersection \(P_C\), so
\[
C\subseteq V_{\Phi}(P_C).
\]
Thus
\[
C=V_{\Phi}(P_C)=\overline{\{P_C\}}.
\]
Therefore \(P_C\) is a generic point of \(C\). The \(T_0\) property implies uniqueness of a generic point. We have proved sobriety and all remaining conditions for a spectral space in the sense of Hochster \cite{Hochster1969}.
\end{proof}

\begin{theorem}\label{thm:order-comparison}
Let \(R\) be a commutative ring with unity, and let \(\rho\) and \(\sigma\) be finite-type nuclei on \(\Id(R)\). If
\[
\rho(I)\subseteq\sigma(I)
\]
for every ideal \(I\lhd R\), then
\[
\Spec_{\sigma}(R)\subseteq\Spec_{\rho}(R),
\]
and the inclusion is a spectral topological embedding.
\end{theorem}

\begin{proof}
Let \(P\in\Spec_{\sigma}(R)\). Then
\[
P\subseteq\rho(P)\subseteq\sigma(P)=P.
\]
Hence \(\rho(P)=P\), so \(P\in\Spec_{\rho}(R)\). This proves the inclusion.

By Proposition~\ref{prop:subspace-realization}, both spaces have the subspace topology inherited from \(\Spec(R)\). Therefore the inclusion is a topological embedding. Its inverse image on a basic quasi-compact open is
\[
D_{\rho}(f)\cap\Spec_{\sigma}(R)=D_{\sigma}(f),
\]
which is quasi-compact by Theorem~\ref{thm:spectrality}. Every quasi-compact open of \(\Spec_{\rho}(R)\) is a finite union of principal opens, because the principal opens form a basis and the given open is quasi-compact. Its inverse image is therefore a finite union of quasi-compact principal opens in \(\Spec_{\sigma}(R)\), hence is quasi-compact. Thus the inclusion is a spectral map.
\end{proof}

\section{Observation Nuclei and Infinitesimal Neighborhoods}\label{sec:observation}

The preceding sections identify the inherited ideal-system structure. We now introduce the point-dependent congruence construction that organizes the affine applications. The general observation nucleus is defined from quotient maps. Its polynomial specialization uses infinitesimal neighborhoods of closed points, and its basic closed-observation case recovers symbolic powers.

\begin{definition}[Observation nucleus]\label{def:observation-nucleus}
Let \(R\) be a commutative ring with unity, and let
\[
\Pi=\{\pi_{\lambda}:R\to Q_{\lambda}\}_{\lambda\in\Lambda}
\]
be a family of unital ring homomorphisms indexed by a set \(\Lambda\). For every subset \(A\subseteq R\), define
\[
\Phiop_{\Pi}(A)
=
\bigcap_{\lambda\in\Lambda}
\pi_{\lambda}^{-1}\bigl(\pi_{\lambda}(A)\bigr).
\]
If \(\Lambda=\varnothing\), the intersection is defined to be \(R\), so that \(\Phiop_{\Pi}(A)=R\) for every \(A\subseteq R\).
\end{definition}

\begin{theorem}\label{thm:observation-validity}
Let \(R\) and \(\Pi\) satisfy Definition~\ref{def:observation-nucleus}. Then \(\Phiop_{\Pi}\) is an algebra-compatible subset closure. Moreover,
\[
\Phiop_{\Pi}(A)
=
\bigcap_{\lambda\in\Lambda}
\bigl(A+\Ker\pi_{\lambda}\bigr)
\]
for every subset \(A\subseteq R\), with the empty intersection interpreted as \(R\).
\end{theorem}

\begin{proof}
If \(\Lambda=\varnothing\), then \(\Phiop_{\Pi}\) is the constant operator with value \(R\). It is extensive, monotone, and idempotent, and both compatibility inclusions have target \(R\). Hence the assertion holds in this case.

Assume that \(\Lambda\neq\varnothing\). Fix \(\lambda\in\Lambda\), and set
\[
K_{\lambda}=\Ker\pi_{\lambda}.
\]
An element \(r\in R\) belongs to \(\pi_{\lambda}^{-1}(\pi_{\lambda}(A))\) if and only if there exists \(a\in A\) such that
\[
\pi_{\lambda}(r)=\pi_{\lambda}(a).
\]
This equality is equivalent to \(r-a\in K_{\lambda}\), and hence to \(r\in A+K_{\lambda}\). Therefore
\[
\pi_{\lambda}^{-1}(\pi_{\lambda}(A))=A+K_{\lambda}.
\]
Intersecting over \(\lambda\) proves the displayed formula.

Extensivity follows from
\[
A\subseteq A+K_{\lambda}
\]
for every \(\lambda\). If \(A\subseteq B\), then
\[
A+K_{\lambda}\subseteq B+K_{\lambda}
\]
for every \(\lambda\), and monotonicity follows after intersection.

We prove idempotence. Extensivity gives
\[
\Phiop_{\Pi}(A)
\subseteq
\Phiop_{\Pi}(\Phiop_{\Pi}(A)).
\]
For each \(\lambda\), the definition gives
\[
\Phiop_{\Pi}(A)\subseteq A+K_{\lambda}.
\]
Consequently,
\[
\Phiop_{\Pi}(A)+K_{\lambda}
\subseteq
A+K_{\lambda}+K_{\lambda}
=
A+K_{\lambda}.
\]
The reverse inclusion follows from \(A\subseteq\Phiop_{\Pi}(A)\). Hence
\[
\Phiop_{\Pi}(A)+K_{\lambda}=A+K_{\lambda}
\]
for every \(\lambda\). Intersecting these equalities yields
\[
\begin{aligned}
\Phiop_{\Pi}(\Phiop_{\Pi}(A))
&=
\bigcap_{\lambda\in\Lambda}
\bigl(\Phiop_{\Pi}(A)+K_{\lambda}\bigr)\\
&=
\bigcap_{\lambda\in\Lambda}(A+K_{\lambda})\\
&=
\Phiop_{\Pi}(A).
\end{aligned}
\]

Let \(x\in\Phiop_{\Pi}(A)\) and \(y\in\Phiop_{\Pi}(B)\). For each \(\lambda\), choose
\[
a_{\lambda}\in A,
\qquad
b_{\lambda}\in B,
\qquad
u_{\lambda},v_{\lambda}\in K_{\lambda}
\]
such that
\[
x=a_{\lambda}+u_{\lambda},
\qquad
y=b_{\lambda}+v_{\lambda}.
\]
Then
\[
x+y=(a_{\lambda}+b_{\lambda})+(u_{\lambda}+v_{\lambda})
\in A+B+K_{\lambda}.
\]
Since this holds for every \(\lambda\),
\[
\Phiop_{\Pi}(A)+\Phiop_{\Pi}(B)
\subseteq
\Phiop_{\Pi}(A+B).
\]
If \(r\in R\), then
\[
rx=ra_{\lambda}+ru_{\lambda}
\in rA+K_{\lambda}
\]
for every \(\lambda\), because \(K_{\lambda}\) is an ideal. Hence
\[
r\Phiop_{\Pi}(A)\subseteq\Phiop_{\Pi}(rA).
\]
All axioms in Definition~\ref{def:algebra-compatible} follow.
\end{proof}

\begin{proposition}\label{prop:nontranslation-criterion}
Let \(R\) be a commutative ring with unity, and let \(K_1,K_2\subsetneq R\) be ideals satisfying \(K_1+K_2=R\). If
\[
\Phiop(A)=(A+K_1)\cap(A+K_2)
\]
for every subset \(A\subseteq R\), then there is no ideal \(J\lhd R\) satisfying
\[
\Phiop(A)=A+J
\]
for every subset \(A\subseteq R\).
\end{proposition}

\begin{proof}
The Chinese remainder theorem \cite[Proposition~1.10]{AtiyahMacdonald1969} gives an element \(e\in R\) such that
\[
e\in K_1,
\qquad
e-1\in K_2.
\]
Therefore
\[
e\in(\{0,1\}+K_1)\cap(\{0,1\}+K_2)=\Phiop(\{0,1\}).
\]
Assume that \(\Phiop(A)=A+J\) for every subset \(A\). Taking \(A=\{0\}\) gives
\[
J=\Phiop(\{0\})=K_1\cap K_2.
\]
The membership \(e\in\{0,1\}+J\) implies either \(e\in J\) or \(e-1\in J\). If \(e\in J\subseteq K_2\), then \(e,e-1\in K_2\), and hence \(1\in K_2\), contrary to \(K_2\subsetneq R\). If \(e-1\in J\subseteq K_1\), then \(e,e-1\in K_1\), and hence \(1\in K_1\), contrary to \(K_1\subsetneq R\). Thus no such ideal \(J\) exists.
\end{proof}

\begin{notation}[Infinitesimal observation algebra]\label{not:jet-algebra}
Let \(k\) be a field, let \(R=k[x_1,\ldots,x_n]\), and let \(a=(a_1,\ldots,a_n)\in k^n\). Define
\[
\m_a=(x_1-a_1,\ldots,x_n-a_n)
\]
and, for every ideal \(I\lhd R\), define
\[
\Vaff(I)=\{a\in k^n:f(a)=0\text{ for every }f\in I\}.
\]
For a positive integer \(q\), define
\[
\mathcal O_{a}^{<q}=R/\m_a^q
\]
and let \(\pi_a^q:R\to\mathcal O_a^{<q}\) be the quotient map.
\end{notation}

\begin{remark}[Infinitesimal neighborhoods and truncated arcs]\label{rem:hasse-jet}
For \(q\geq1\), the Artin local algebra \(R/\m_a^q\) is the coordinate algebra of the \((q-1)\)-st infinitesimal neighborhood of \(a\) in affine \(n\)-space. Equivalently, the residue class of a polynomial is determined by its expansion in the shifted variables \(x_i-a_i\) through total degree less than \(q\); these coefficients may be expressed by Hasse derivatives in arbitrary characteristic. This algebra is not the coordinate algebra representing the scheme of truncated arcs. Hasse--Schmidt algebras and scheme-theoretic jets are treated by Vojta \cite{Vojta2007}; see also Musta\c{t}\u{a} \cite{Mustata2001}.

Local jet closures constitute a distinct construction. De Fernex, Ein, and Ishii define them by equality of local jet schemes or jet supports \cite{deFernexEinIshii2018}, and Chen and Zuo study their associated algebras, filtrations, and special homogeneous cases \cite{ChenZuo2025}. The factor \(I+\m_a^q\) used here is an elementary saturation by an infinitesimal-neighborhood congruence. No general equality with a local jet closure is asserted.

There is an exact coincidence under additional hypotheses. Let $k$ be algebraically closed of characteristic zero, let $B=k[[z_1,\ldots,z_n]]$ with maximal ideal $\mathfrak n=(z_1,\ldots,z_n)$, and let $J$ be extended from a homogeneous polynomial ideal. Chen and Zuo \cite[Theorem~A]{ChenZuo2025} prove that its $m$-th local jet closure satisfies
\[
J^{m\text{-}\mathrm{jc}}=J+\mathfrak n^{m+1}.
\]
Thus the completed congruence factor agrees with their jet closure in this case. For the limiting arc closure, Mallory \cite{Mallory2020} proves triviality on Noetherian local $k$-algebras whose residue-field extension is separable. These local reconstruction results complement the global closed-point recovery theorem used in Section~\ref{sec:tower}.
\end{remark}

\begin{definition}[Point-dependent infinitesimal observation]\label{def:jet-observation}
Let \(k\) be a field, and let \(R=k[x_1,\ldots,x_n]\). Let \(X\subseteq k^n\), and let \(\nu:X\to\bbN_{>0}\). For every subset \(A\subseteq R\), define
\[
\Phiop_{X,\nu}(A)
=
\bigcap_{a\in X}
(\pi_a^{\nu(a)})^{-1}
\bigl(\pi_a^{\nu(a)}(A)\bigr)
=
\bigcap_{a\in X}
\bigl(A+\m_a^{\nu(a)}\bigr).
\]
If \(X=\varnothing\), the intersection is defined to be \(R\). If \(\nu(a)=q\) for every \(a\in X\), write \(\Phiop_{X,q}\). The function \(\nu\) is bounded if there exists a positive integer \(N\) such that \(\nu(a)\leq N\) for every \(a\in X\).
\end{definition}

\begin{corollary}\label{cor:jet-validity}
Let \(k\), \(R\), \(X\), and \(\nu\) satisfy Definition~\ref{def:jet-observation}. Then \(\Phiop_{X,\nu}\) is an algebra-compatible subset closure. Its restriction to \(\Id(R)\) is a finite-type nucleus. If \(X=\varnothing\), its fixed-prime space is empty.
\end{corollary}

\begin{proof}
If \(X\neq\varnothing\), Theorem~\ref{thm:observation-validity} applies to the quotient maps \(\pi_a^{\nu(a)}\). If \(X=\varnothing\), the operator is constant with value \(R\), so the axioms follow directly. In either case, the Hilbert basis theorem implies that \(R\) is Noetherian \cite{Eisenbud1995}. Every ideal is therefore finitely generated, and the restricted nucleus is of finite type. In the empty case, every ideal \(P\) satisfies \(\Phiop_{X,\nu}(P)=R\), so no proper prime ideal is fixed.
\end{proof}

\begin{proposition}\label{prop:symbolic-powers}
Let \(k\) be algebraically closed, let \(R=k[x_1,\ldots,x_n]\), let \(P\in\Spec(R)\), and let \(q\geq1\). Then
\[
\Phiop_{\Vaff(P),q}(0)
=
\bigcap_{a\in\Vaff(P)}\m_a^q
=
P^{(q)},
\]
where \(P^{(q)}=P^qR_P\cap R\) is the \(q\)-th symbolic power of \(P\).
\end{proposition}

\begin{proof}
The maximal ideals containing \(P\) are precisely the point ideals \(\m_a\) with \(a\in\Vaff(P)\), by Hilbert's Nullstellensatz \cite{Eisenbud1995}. The defining formula for \(\Phiop_{\Vaff(P),q}(0)\) gives the first equality. Since the algebraically closed field \(k\) is perfect, the Zariski--Nagata theorem identifies the intersection of the \(q\)-th powers of all closed-point maximal ideals containing \(P\) with the symbolic power \(P^{(q)}\); see \cite{DaoEtAl2018}. This proves the second equality.
\end{proof}

\section{The Exact Relative Nullstellensatz}\label{sec:nullstellensatz}

The observation nucleus records both the selected point set and the order assigned to each point. This section separates those two parameters. A bounded order function supplies a single exponent valid at all observed zeros, which determines the radical exactly. Algebraic closedness is not needed for this uniform-exponent argument; it enters only when affine point sets are identified with radical ideals and when generic points are reconstructed.

\begin{notation}[Relative vanishing data]\label{not:relative-vanishing}
Let \(k\) be a field, let \(R=k[x_1,\ldots,x_n]\), and let \(X\subseteq k^n\). For an ideal \(I\lhd R\), define
\[
V_X(I)=X\cap\Vaff(I).
\]
For a subset \(Y\subseteq k^n\), define
\[
\calI(Y)=\{f\in R:f(a)=0\text{ for every }a\in Y\}.
\]
In particular, \(\calI(\varnothing)=R\).
\end{notation}

\begin{theorem}\label{thm:relative-nullstellensatz}
Let \(k\) be a field, let \(R=k[x_1,\ldots,x_n]\), let \(X\subseteq k^n\), and let \(\nu:X\to\bbN_{>0}\) be bounded. Then, for every ideal \(I\lhd R\),
\[
\sqrt{\Phiop_{X,\nu}(I)}
=
\calI\bigl(V_X(I)\bigr).
\]
\end{theorem}

\begin{proof}
Choose a positive integer \(N\) such that
\[
\nu(a)\leq N
\]
for every \(a\in X\). If \(X=\varnothing\), then \(\Phiop_{X,\nu}(I)=R\) and \(V_X(I)=\varnothing\). Hence both sides of the asserted equality are \(R\). Assume from now on that \(X\neq\varnothing\).

Let \(a\in X\setminus\Vaff(I)\). There exists \(g\in I\) such that \(g(a)\neq0\). Set
\[
c=g(a)\in k^{\times},
\qquad
h=g-c.
\]
Since \(h(a)=0\), we have \(h\in\m_a\). Put \(q=\nu(a)\). In \(R/\m_a^q\), the class of \(h\) is nilpotent and satisfies \(h^q=0\). The finite geometric sum
\[
c^{-1}\sum_{j=0}^{q-1}(-c^{-1}h)^j
\]
is therefore an inverse for the class of \(c+h=g\). Since \(g\in I\), the image of \(I\) in \(R/\m_a^q\) contains a unit. Consequently,
\[
I+\m_a^{\nu(a)}=R.
\tag{5.1}\label{eq:outside-factor}
\]

Let \(a\in V_X(I)\). Every element of \(I\) vanishes at \(a\), so \(I\subseteq\m_a\). Since \(\nu(a)\geq1\),
\[
\m_a^{\nu(a)}\subseteq\m_a.
\]
Hence
\[
I\subseteq I+\m_a^{\nu(a)}\subseteq\m_a.
\tag{5.2}\label{eq:inside-sandwich}
\]
Equation~\eqref{eq:outside-factor} shows that every point outside \(V_X(I)\) contributes the unit ideal. Thus
\[
\Phiop_{X,\nu}(I)
=
\bigcap_{a\in V_X(I)}
\bigl(I+\m_a^{\nu(a)}\bigr),
\tag{5.3}\label{eq:reduced-intersection}
\]
where the intersection is \(R\) when \(V_X(I)=\varnothing\).

Let
\[
f\in\sqrt{\Phiop_{X,\nu}(I)}.
\]
There exists a positive integer \(M\) such that
\[
f^M\in\Phiop_{X,\nu}(I).
\]
Fix \(a\in V_X(I)\). Equations~\eqref{eq:reduced-intersection} and \eqref{eq:inside-sandwich} give
\[
f^M\in I+\m_a^{\nu(a)}\subseteq\m_a.
\]
The maximal ideal \(\m_a\) is prime, so \(f\in\m_a\). Thus \(f(a)=0\). Since \(a\) was arbitrary,
\[
f\in\calI(V_X(I)).
\]
Therefore
\[
\sqrt{\Phiop_{X,\nu}(I)}
\subseteq
\calI(V_X(I)).
\]

Conversely, let \(f\in\calI(V_X(I))\). For every \(a\in V_X(I)\), one has \(f\in\m_a\), and hence
\[
f^N\in\m_a^N.
\]
The inequality \(N\geq\nu(a)\) gives
\[
\m_a^N\subseteq\m_a^{\nu(a)}.
\]
Consequently,
\[
f^N\in\m_a^{\nu(a)}
\subseteq
I+\m_a^{\nu(a)}
\]
for every \(a\in V_X(I)\). Equation~\eqref{eq:reduced-intersection} yields
\[
f^N\in\Phiop_{X,\nu}(I).
\]
Thus
\[
f\in\sqrt{\Phiop_{X,\nu}(I)}.
\]
The reverse inclusion follows, and the proof is complete.
\end{proof}

\begin{example}[Varying observation order along a parabola]\label{ex:parabola-orders}
Let $k$ be algebraically closed, let $R=k[x,y]$, and put
\[
h=y-x^2,\qquad
X=\Vaff((h))=\{a_t=(t,t^2):t\in k\},\qquad
I=(h^3(h-1)).
\]
Assign the bounded, nonconstant observation order
\[
\nu(a_t)=
\begin{cases}
3,&t=0,\\
2,&t\neq0.
\end{cases}
\]
We compute the observation ideal, compare it with constant orders, and distinguish its radical from both $I$ and $\sqrt I$.

The change of coordinates $h=y-x^2$ identifies $R$ with $k[x,h]$. Since
\[
y-t^2=h+(x-t)(x+t),
\]
we have $\m_{a_t}=(x-t,h)$ in every characteristic. Moreover,
\[
I\subseteq(h^3)\subseteq\m_{a_t}^{\nu(a_t)},
\]
so every factor $I+\m_{a_t}^{\nu(a_t)}$ equals $\m_{a_t}^{\nu(a_t)}$.

For every positive integer $q$ and either $S=k$ or $S=k\setminus\{0\}$, one has
\[
\bigcap_{t\in S}(x-t,h)^q=(h^q).
\]
To prove the nontrivial inclusion, write a polynomial in the intersection uniquely as
\[
f=\sum_{j=0}^{q-1}b_j(x)h^j+h^q g(x,h),
\qquad b_j(x)\in k[x].
\]
For each $t\in S$, the homomorphism $k[x,h]\to k[h]/(h^q)$ sending $x$ to $t$ sends $(x-t,h)^q$ to zero. Hence
\[
\sum_{j=0}^{q-1}b_j(t)h^j=0\quad\text{in }k[h]/(h^q).
\]
The residue classes of $1,h,\ldots,h^{q-1}$ are linearly independent, so $b_j(t)=0$ for every $j<q$ and every $t\in S$. An algebraically closed field is infinite, and both choices of $S$ are infinite. Therefore each $b_j$ is the zero polynomial, and $f\in(h^q)$. The reverse inclusion follows from $h\in(x-t,h)$.

Applying this identity away from the origin gives
\[
\begin{aligned}
\Phiop_{X,\nu}(I)
&=(h^2)\cap(x,h)^3\\
&=h^2(x,h)=(xh^2,h^3)\\
&=\bigl(x(y-x^2)^2,(y-x^2)^3\bigr).
\end{aligned}
\]
For the middle equality, write $g=c+xu+hv$, with $c\in k$ and $u,v\in k[x,h]$. Then $h^2g$ lies in $(x,h)^3$ if and only if $c=0$, since its only possible term of total degree less than three in the coordinates $x,h$ is $ch^2$. This proves $(h^2)\cap(x,h)^3=h^2(x,h)$.

For the constant observation orders two and three, the same intersection identity yields
\[
\Phiop_{X,2}(I)=(h^2),\qquad
\Phiop_{X,3}(I)=(h^3).
\]
Thus increasing the order from two to three only at the origin changes $(h^2)$ into $(xh^2,h^3)$. The strict containments
\[
I\subsetneq(h^3)\subsetneq(xh^2,h^3)\subsetneq(h^2)\subsetneq(h)
\]
are witnessed, respectively, by $h^3$, $xh^2$, $h^2$, and $h$.

On the other hand, $V_X(I)=X$, and the intersection identity for $q=1$ gives $\calI(X)=(h)$. Since $(h^3)\subseteq\Phiop_{X,\nu}(I)\subseteq(h)$ and $(h)$ is prime, we obtain
\[
\sqrt{\Phiop_{X,\nu}(I)}
=\calI(V_X(I))=(h)\neq I.
\]
In fact, the coprime prime factors $h$ and $h-1$ give
\[
\sqrt I=(h(h-1))\subsetneq(h).
\]
The observation set omits the component $\Vaff((h-1))$ of $\Vaff(I)$. The nonconstant order changes the retained infinitesimal conditions on the remaining parabola, while $\calI(V_X(I))$ is independent of $\nu$.
\end{example}

\begin{remark}\label{rem:field-independent-nullstellensatz}
Theorem~\ref{thm:relative-nullstellensatz} is valid over an arbitrary field because it concerns only selected \(k\)-rational points and powers of their maximal ideals. Algebraic closedness is required in the following results when Hilbert's Nullstellensatz \cite{Eisenbud1995} is used to recover ideals and ambient irreducible varieties from their point sets.
\end{remark}

\begin{theorem}\label{thm:relative-correspondence}
Let \(k\) be a field, let \(R=k[x_1,\ldots,x_n]\), let \(X\subseteq k^n\), and let \(\nu:X\to\bbN_{>0}\) be bounded. Define
\[
\mathfrak r_X(I)=\sqrt{\Phiop_{X,\nu}(I)}
\]
for every ideal \(I\lhd R\). Then \(\mathfrak r_X\) is an extensive, monotone, idempotent radical closure operation on \(\Id(R)\), it is independent of \(\nu\), and the assignments
\[
I\longmapsto V_X(I)
\qquad\text{and}\qquad
Y\longmapsto\calI(Y)
\]
induce mutually inverse order-reversing bijections between the \(\mathfrak r_X\)-closed ideals and the Zariski closed subsets of the subspace \(X\).
\end{theorem}

\begin{proof}
Theorem~\ref{thm:relative-nullstellensatz} gives
\[
\mathfrak r_X(I)=\calI(V_X(I)).
\tag{5.4}\label{eq:relative-radicalized}
\]
The right-hand side is independent of \(\nu\), and it is a radical ideal because every vanishing ideal is radical.

If \(f\in I\), then \(f\) vanishes at every point of \(V_X(I)\). Hence
\[
I\subseteq\calI(V_X(I))=\mathfrak r_X(I),
\]
proving extensivity. If \(I\subseteq J\), then
\[
V_X(J)\subseteq V_X(I).
\]
Taking vanishing ideals reverses inclusion, so
\[
\mathfrak r_X(I)=\calI(V_X(I))
\subseteq
\calI(V_X(J))=\mathfrak r_X(J).
\]
Thus \(\mathfrak r_X\) is monotone.

Extensivity gives
\[
V_X(\mathfrak r_X(I))\subseteq V_X(I).
\]
Conversely, every element of \(\mathfrak r_X(I)=\calI(V_X(I))\) vanishes on \(V_X(I)\), and therefore
\[
V_X(I)\subseteq V_X(\mathfrak r_X(I)).
\]
Thus
\[
V_X(\mathfrak r_X(I))=V_X(I).
\]
Applying \(\calI\) and using \eqref{eq:relative-radicalized} gives
\[
\mathfrak r_X(\mathfrak r_X(I))
=
\calI(V_X(\mathfrak r_X(I)))
=
\calI(V_X(I))
=
\mathfrak r_X(I).
\]
Hence \(\mathfrak r_X\) is idempotent.

Let \(Y\) be closed in the subspace \(X\). There exists an ideal \(I\lhd R\) such that
\[
Y=X\cap\Vaff(I)=V_X(I).
\]
By definition, every polynomial in \(\calI(Y)\) vanishes on \(Y\), so
\[
Y\subseteq V_X(\calI(Y)).
\]
Since \(I\subseteq\calI(Y)\), one also has
\[
V_X(\calI(Y))\subseteq V_X(I)=Y.
\]
Therefore
\[
V_X(\calI(Y))=Y.
\tag{5.5}\label{eq:relative-set-recovery}
\]
If \(J\) is \(\mathfrak r_X\)-closed, then
\[
\calI(V_X(J))=\mathfrak r_X(J)=J.
\tag{5.6}\label{eq:relative-ideal-recovery}
\]
Equations~\eqref{eq:relative-set-recovery} and \eqref{eq:relative-ideal-recovery} give mutually inverse correspondences. Their order-reversing character follows from the definitions.
\end{proof}

\begin{corollary}\label{cor:full-nullstellensatz}
Let \(k\) be algebraically closed, let \(R=k[x_1,\ldots,x_n]\), and let \(\nu:k^n\to\bbN_{>0}\) be bounded. Then, for every ideal \(I\lhd R\),
\[
\sqrt{\Phiop_{k^n,\nu}(I)}
=
\calI(\Vaff(I))
=
\sqrt I.
\]
\end{corollary}

\begin{proof}
The first equality is Theorem~\ref{thm:relative-nullstellensatz} with \(X=k^n\). The second equality is Hilbert's Nullstellensatz \cite{Eisenbud1995}.
\end{proof}

\begin{theorem}\label{thm:density-characterization}
Let \(k\) be algebraically closed, let \(R=k[x_1,\ldots,x_n]\), let \(X\subseteq k^n\), and let \(\nu:X\to\bbN_{>0}\) be bounded. Then, for every prime ideal \(P\in\Spec(R)\), the following conditions are equivalent:
\[
\Phiop_{X,\nu}(P)=P,
\]
\[
\calI\bigl(X\cap\Vaff(P)\bigr)=P,
\]
and
\[
\overline{X\cap\Vaff(P)}^{\Zar}=\Vaff(P).
\]
Consequently,
\[
\Spec_{\Phi_{X,\nu}}(R)
=
\{P\in\Spec(R):X\cap\Vaff(P)\text{ is Zariski dense in }\Vaff(P)\}.
\]
\end{theorem}

\begin{proof}
Assume that \(\Phiop_{X,\nu}(P)=P\). Since \(P\) is radical, Theorem~\ref{thm:relative-nullstellensatz} gives
\[
P
=
\sqrt P
=
\sqrt{\Phiop_{X,\nu}(P)}
=
\calI(X\cap\Vaff(P)).
\]
Thus the first condition implies the second.

Assume that
\[
\calI(X\cap\Vaff(P))=P.
\]
For \(a\in X\setminus\Vaff(P)\), equation~\eqref{eq:outside-factor} gives
\[
P+\m_a^{\nu(a)}=R.
\]
For \(a\in X\cap\Vaff(P)\), equation~\eqref{eq:inside-sandwich} gives
\[
P\subseteq P+\m_a^{\nu(a)}\subseteq\m_a.
\]
Intersecting yields
\[
P
\subseteq
\Phiop_{X,\nu}(P)
\subseteq
\bigcap_{a\in X\cap\Vaff(P)}\m_a
=
\calI(X\cap\Vaff(P))
=P.
\]
Hence \(\Phiop_{X,\nu}(P)=P\), proving the equivalence of the first two conditions.

For every subset \(Y\subseteq k^n\), Hilbert's Nullstellensatz \cite{Eisenbud1995} gives
\[
\overline{Y}^{\Zar}=\Vaff(\calI(Y))
\qquad\text{and}\qquad
\calI(Y)=\calI(\overline{Y}^{\Zar}).
\]
Take \(Y=X\cap\Vaff(P)\). If \(\calI(Y)=P\), then
\[
\overline{Y}^{\Zar}=\Vaff(\calI(Y))=\Vaff(P).
\]
Conversely, if \(\overline{Y}^{\Zar}=\Vaff(P)\), then Hilbert's Nullstellensatz \cite{Eisenbud1995} gives
\[
\calI(Y)
=
\calI(\overline{Y}^{\Zar})
=
\calI(\Vaff(P))
=
\sqrt P
=
P.
\]
This proves the equivalence with the density condition and the displayed description of the fixed-prime spectrum.
\end{proof}

\begin{definition}[Sobrification]\label{def:sobrification}
Let \(T\) be a topological space. Define \(\Sob(T)\) to be the set of nonempty irreducible closed subsets of \(T\). For every open subset \(U\subseteq T\), define
\[
\widetilde U=\{C\in\Sob(T):C\cap U\neq\varnothing\}.
\]
Equip \(\Sob(T)\) with the topology generated by the sets \(\widetilde U\).
\end{definition}

\begin{theorem}\label{thm:sobrification}
Let \(k\) be algebraically closed, and let \(R=k[x_1,\ldots,x_n]\). Let \(X\subseteq k^n\) carry the induced Zariski topology, and let \(\nu:X\to\bbN_{>0}\) be bounded. Then the map
\[
\Theta:\Spec_{\Phi_{X,\nu}}(R)\longrightarrow\Sob(X),
\qquad
\Theta(P)=X\cap\Vaff(P),
\]
is a homeomorphism.
\end{theorem}

\begin{proof}
Let \(P\in\Spec_{\Phi_{X,\nu}}(R)\). The affine variety \(\Vaff(P)\) is irreducible because \(P\) is prime \cite{Eisenbud1995}. By Theorem~\ref{thm:density-characterization}, the subset
\[
X\cap\Vaff(P)
\]
is dense in \(\Vaff(P)\), and in particular it is nonempty. A dense subspace of an irreducible space is irreducible \cite{GortzWedhorn2020}. Moreover, \(X\cap\Vaff(P)\) is closed in \(X\). Hence \(\Theta(P)\) is a nonempty irreducible closed subset of \(X\), so \(\Theta\) is well defined.

If \(P,Q\in\Spec_{\Phi_{X,\nu}}(R)\) and
\[
X\cap\Vaff(P)=X\cap\Vaff(Q),
\]
then Theorem~\ref{thm:density-characterization} gives
\[
P=\calI(X\cap\Vaff(P))
=\calI(X\cap\Vaff(Q))=Q.
\]
Thus \(\Theta\) is injective.

Let \(C\in\Sob(X)\). Since \(C\) is closed in \(X\), there exists a Zariski closed subset \(F\subseteq k^n\) such that
\[
C=X\cap F.
\]
Let \(\overline C\) be the Zariski closure of \(C\) in \(k^n\). The closure of an irreducible subset is irreducible \cite{GortzWedhorn2020}, so
\[
\overline C=\Vaff(P)
\]
for the prime ideal \(P=\calI(C)\). Since \(C\subseteq F\), one has \(\overline C\subseteq F\), and therefore
\[
X\cap\Vaff(P)
\subseteq
X\cap F=C.
\]
The reverse inclusion follows from \(C\subseteq\overline C=\Vaff(P)\). Hence
\[
C=X\cap\Vaff(P).
\]
The set \(C\) is dense in \(\Vaff(P)\) by construction. Theorem~\ref{thm:density-characterization} gives
\[
P\in\Spec_{\Phi_{X,\nu}}(R),
\]
so \(\Theta(P)=C\). Thus \(\Theta\) is surjective.

For \(f\in R\), let
\[
U_f=X\cap D(f),
\qquad
D(f)=\{a\in k^n:f(a)\neq0\}.
\]
The sets \(U_f\) form a basis of the induced Zariski topology on \(X\). We claim that
\[
\Theta^{-1}(\widetilde{U_f})=D_{\Phi}(f).
\tag{5.7}\label{eq:sob-basic}
\]
Let \(P\in\Spec_{\Phi_{X,\nu}}(R)\). Suppose first that \(P\in D_{\Phi}(f)\), so \(f\notin P\). The open subset \(D(f)\cap\Vaff(P)\) is nonempty. Indeed, if it were empty, then \(f\) would vanish on \(\Vaff(P)\), and Hilbert's Nullstellensatz \cite{Eisenbud1995} would imply
\[
f\in\calI(\Vaff(P))=P,
\]
contrary to \(f\notin P\). Since \(X\cap\Vaff(P)\) is dense in \(\Vaff(P)\), it meets the nonempty open subset \(D(f)\cap\Vaff(P)\). Therefore
\[
(X\cap\Vaff(P))\cap U_f\neq\varnothing,
\]
and \(\Theta(P)\in\widetilde{U_f}\).

Conversely, suppose that \(\Theta(P)\cap U_f\neq\varnothing\). There exists \(a\in X\cap\Vaff(P)\) such that \(f(a)\neq0\). Hence \(f\notin\m_a\). Since \(P\subseteq\m_a\), it follows that \(f\notin P\), and therefore \(P\in D_{\Phi}(f)\). This proves \eqref{eq:sob-basic}.

The sets \(D_{\Phi}(f)\) form a basis of the source, and the sets \(\widetilde{U_f}\) form a basis of the target. Equation~\eqref{eq:sob-basic}, together with bijectivity, proves that \(\Theta\) and its inverse are continuous. Hence \(\Theta\) is a homeomorphism.
\end{proof}

\begin{remark}\label{rem:sobrification-context}
Theorem~\ref{thm:sobrification} realizes a general topological construction. For a subspace $S$ of a sober space $T$, its sobrification is homeomorphic to
\[
\left\{t\in T:
\overline{S\cap\overline{\{t\}}}^{\,T}
=\overline{\{t\}}^{\,T}\right\}
\]
with the subspace topology; this is the Goldman closure of $S$ \cite[Proposition~4.1 and Theorem~4.20]{BouacidaEtAl2003}. See also \cite[Sections~2.2--2.3]{Schwartz2023}. Taking $T=\Spec(R)$ and $S=\{\m_a:a\in X\}$ identifies this locus with the fixed primes through Theorem~\ref{thm:density-characterization}. Recovery from all closed points is the familiar special case \cite{GortzWedhorn2020}.
\end{remark}

\begin{corollary}\label{cor:jet-spectrality}
Let \(k\) be algebraically closed, let \(R=k[x_1,\ldots,x_n]\), let \(X\subseteq k^n\), and let \(\nu:X\to\bbN_{>0}\) be bounded. Then \(\Spec_{\Phi_{X,\nu}}(R)\) is a spectral space.
\end{corollary}

\begin{proof}
Corollary~\ref{cor:jet-validity} shows that the restricted nucleus is of finite type. The assertion follows from Theorem~\ref{thm:spectrality}.
\end{proof}

\begin{corollary}\label{cor:observation-loci}
Let \(k\) be algebraically closed, let \(R=k[x_1,\ldots,x_n]\), and let \(\nu\) be bounded on the indicated observation set. Then the following conclusions hold.

If \(X=k^n\), then
\[
\Spec_{\Phi_{X,\nu}}(R)=\Spec(R).
\]
If \(X=\Vaff(J)\) for an ideal \(J\lhd R\), then
\[
\Spec_{\Phi_{X,\nu}}(R)
=
\VSpec(J)
:=
\{P\in\Spec(R):J\subseteq P\}.
\]
If \(X\subseteq k^n\) is finite, then
\[
\Spec_{\Phi_{X,\nu}}(R)=\{\m_a:a\in X\}.
\]
\end{corollary}

\begin{proof}
If \(X=k^n\), then
\[
X\cap\Vaff(P)=\Vaff(P)
\]
for every prime \(P\), and Theorem~\ref{thm:density-characterization} gives the first equality.

Let \(X=\Vaff(J)\). For a prime \(P\),
\[
X\cap\Vaff(P)=\Vaff(J+P).
\]
If \(J\subseteq P\), this set equals \(\Vaff(P)\) and is dense. Conversely, if it is dense in \(\Vaff(P)\), then Theorem~\ref{thm:density-characterization} and Hilbert's Nullstellensatz \cite{Eisenbud1995} give
\[
P
=
\calI(\Vaff(J+P))
=
\sqrt{J+P},
\]
which implies \(J\subseteq P\). Hence the fixed primes are precisely the elements of \(\VSpec(J)\).

Assume that \(X\) is finite. If \(a\in X\), then
\[
X\cap\Vaff(\m_a)=\{a\}=\Vaff(\m_a),
\]
so \(\m_a\) is fixed. Conversely, let \(P\) be fixed. The finite set \(X\cap\Vaff(P)\) is dense in the irreducible closed set \(\Vaff(P)\). Every point of \(k^n\) is Zariski closed by Hilbert's Nullstellensatz \cite{Eisenbud1995}, so a finite subset is closed. Therefore
\[
\Vaff(P)=X\cap\Vaff(P).
\]
A finite irreducible space whose points are closed has exactly one point: otherwise one point and the finite union of the remaining point closures would be two proper closed subsets whose union is the whole space. Thus \(\Vaff(P)=\{a\}\) for some \(a\in X\), and Hilbert's Nullstellensatz \cite{Eisenbud1995} gives \(P=\m_a\).
\end{proof}

\begin{theorem}\label{thm:variation}
Let \(k\) be algebraically closed, let \(R=k[x_1,\ldots,x_n]\), and let \(X\subseteq Y\subseteq k^n\). Let \(\nu:X\to\bbN_{>0}\) and \(\mu:Y\to\bbN_{>0}\) be bounded.

If \(\mu|_X=\nu\), then
\[
\Phiop_{Y,\mu}(I)\subseteq\Phiop_{X,\nu}(I)
\]
for every ideal \(I\), and
\[
\Spec_{\Phi_{X,\nu}}(R)
\subseteq
\Spec_{\Phi_{Y,\mu}}(R).
\]

If \(X=Y\) and \(\nu(a)\leq\mu(a)\) for every \(a\in X\), then
\[
\Phiop_{X,\mu}(I)\subseteq\Phiop_{X,\nu}(I)
\]
for every ideal \(I\), while
\[
\Spec_{\Phi_{X,\mu}}(R)
=
\Spec_{\Phi_{X,\nu}}(R).
\]
\end{theorem}

\begin{proof}
Assume first that \(X\subseteq Y\) and \(\mu|_X=\nu\). The intersection defining \(\Phiop_{Y,\mu}(I)\) contains every factor in the intersection defining \(\Phiop_{X,\nu}(I)\), together with the factors indexed by \(Y\setminus X\). Hence
\[
\Phiop_{Y,\mu}(I)\subseteq\Phiop_{X,\nu}(I).
\]
If \(P\) is fixed by \(\Phiop_{X,\nu}\), then \(X\cap\Vaff(P)\) is dense in \(\Vaff(P)\). Since
\[
X\cap\Vaff(P)\subseteq Y\cap\Vaff(P),
\]
the latter set is also dense. Theorem~\ref{thm:density-characterization} gives the inclusion of fixed-prime spectra.

Assume next that \(X=Y\) and \(\nu(a)\leq\mu(a)\) for every \(a\in X\). Then
\[
\m_a^{\mu(a)}\subseteq\m_a^{\nu(a)}
\]
for every \(a\), and therefore
\[
I+\m_a^{\mu(a)}
\subseteq
I+\m_a^{\nu(a)}.
\]
Intersecting gives the inclusion of closures. The density characterization depends only on \(X\), not on the bounded order function, so the fixed-prime spectra are equal.
\end{proof}

\section{Symbolic-Power Interfaces and Infinitesimal Recovery}\label{sec:tower}

The full affine observation model retains every prime ideal while changing nonreduced structure. We first express the quotient of an observation closure as an intersection of powers of closed-point maximal ideals. The resulting tower is governed by the classical nilpotent Nullstellensatz of Eisenbud and Hochster \cite{EisenbudHochster1979}: it stabilizes at the original ideal after finitely many stages. We formulate this theorem in observation language, retain the Krull-intersection argument as an independent local comparison, strengthen nontranslation to every positive affine dimension, and distinguish the finite fat-point case from the infinite observation case.

\begin{definition}[Full affine observation tower]\label{def:jet-tower}
Let \(k\) be algebraically closed, and let \(R=k[x_1,\ldots,x_n]\). Let \(q\geq1\). For every ideal \(I\lhd R\), define
\[
\Phiop_q(I)
=
\bigcap_{a\in k^n}(I+\m_a^q).
\]
\end{definition}

\begin{proposition}\label{prop:quotient-power-interpretation}
Let \(k\) be algebraically closed, let \(R=k[x_1,\ldots,x_n]\), let \(I\lhd R\), let \(A=R/I\), and let \(q\geq1\). Then
\[
\frac{\Phiop_q(I)}{I}
=
\bigcap_{\mathfrak n\in\operatorname{Max}(A)}\mathfrak n^q.
\]
\end{proposition}

\begin{proof}
Equation~\eqref{eq:outside-factor} gives
\[
\Phiop_q(I)
=
\bigcap_{a\in\Vaff(I)}(I+\m_a^q),
\]
because every point outside \(\Vaff(I)\) contributes the unit ideal. All ideals in this intersection contain \(I\), so passage to the quotient commutes with the intersection:
\[
\frac{\Phiop_q(I)}{I}
=
\bigcap_{a\in\Vaff(I)}\frac{I+\m_a^q}{I}.
\]
For \(a\in\Vaff(I)\), the ideal \(\m_a/I\) is a maximal ideal of \(A\), and
\[
\frac{I+\m_a^q}{I}
=
\left(\frac{\m_a}{I}\right)^q.
\]
Conversely, the weak Nullstellensatz \cite{Eisenbud1995} identifies every maximal ideal of the finitely generated \(k\)-algebra \(A\) with \(\m_a/I\) for a unique point \(a\in\Vaff(I)\). Substitution gives the asserted formula.
\end{proof}

\begin{theorem}\label{thm:finite-stage-recovery}
Let \(k\) be algebraically closed, let \(R=k[x_1,\ldots,x_n]\), and let \(I\lhd R\). Then, for every positive integer \(q\),
\[
I\subseteq\Phiop_{q+1}(I)
\subseteq\Phiop_q(I)
\subseteq\sqrt I,
\]
\[
\sqrt{\Phiop_q(I)}=\sqrt I,
\qquad
\Phiop_1(I)=\sqrt I.
\]
Moreover, there exists an integer \(e=e(I)\geq1\) such that
\[
\Phiop_q(I)=I
\qquad
\text{for every }q\geq e.
\]
Consequently,
\[
\bigcap_{q\geq1}\Phiop_q(I)=I.
\]
\end{theorem}

\begin{proof}
Extensivity gives
\[
I\subseteq\Phiop_{q+1}(I).
\]
For every \(a\in k^n\),
\[
\m_a^{q+1}\subseteq\m_a^q,
\]
and hence
\[
I+\m_a^{q+1}\subseteq I+\m_a^q.
\]
Intersecting over all points gives
\[
\Phiop_{q+1}(I)\subseteq\Phiop_q(I).
\]
Corollary~\ref{cor:full-nullstellensatz} gives
\[
\sqrt{\Phiop_q(I)}=\sqrt I.
\]
Since every ideal is contained in its radical,
\[
\Phiop_q(I)
\subseteq
\sqrt{\Phiop_q(I)}
=
\sqrt I.
\]

For \(q=1\), every point outside \(\Vaff(I)\) contributes \(R\), while every point \(a\in\Vaff(I)\) contributes
\[
I+\m_a=\m_a.
\]
Therefore
\[
\Phiop_1(I)
=
\bigcap_{a\in\Vaff(I)}\m_a
=
\calI(\Vaff(I))
=
\sqrt I
\]
by Hilbert's Nullstellensatz \cite{Eisenbud1995}.

For every positive integer \(q\), equation~\eqref{eq:outside-factor} gives
\[
\Phiop_q(I)
=
\bigcap_{a\in\Vaff(I)}(I+\m_a^q).
\tag{6.1}\label{eq:full-affine-zero-intersection}
\]
By the Nullstellensatz with nilpotents of Eisenbud and Hochster \cite[Corollary~2]{EisenbudHochster1979}, there exists an integer \(e=e(I)\geq1\) such that
\[
I
=
\bigcap_{a\in\Vaff(I)}(I+\m_a^e)
=
\Phiop_e(I).
\]
If \(q\geq e\), then \(\m_a^q\subseteq\m_a^e\) for every \(a\in k^n\). Hence
\[
I
\subseteq
\Phiop_q(I)
\subseteq
\Phiop_e(I)
=I.
\]
Thus \(\Phiop_q(I)=I\) for every \(q\geq e\), and the intersection identity follows.

For comparison, the intersection identity also has an independent local proof from the Krull intersection theorem. Assume that \(I\subsetneq R\), let \(f\in R\setminus I\), and set
\[
A=R/I,
\qquad
\overline f=f+I.
\]
The annihilator \(\Ann_A(\overline f)\) is proper. Choose a maximal ideal \(\mathfrak n\) containing it. Then \(\overline f/1\neq0\) in \(A_{\mathfrak n}\), because an equality \(\overline f/1=0\) would provide an element of \(A\setminus\mathfrak n\) annihilating \(\overline f\), contrary to the choice of \(\mathfrak n\). The weak Nullstellensatz \cite{Eisenbud1995} gives
\[
\mathfrak n=\m_a/I
\]
for some \(a\in\Vaff(I)\). The localization \(A_{\mathfrak n}\) is a Noetherian local ring \cite{AtiyahMacdonald1969} and satisfies
\[
\bigcap_{q\geq1}(\mathfrak nA_{\mathfrak n})^q=0
\]
by the Krull intersection theorem \cite[Theorem~10.17]{AtiyahMacdonald1969}. Therefore
\[
\overline f/1\notin(\mathfrak nA_{\mathfrak n})^q
\]
for some positive integer \(q\). If \(f\in I+\m_a^q\), then
\[
\overline f\in(\m_a^q+I)/I=\mathfrak n^q,
\]
and localization would give
\[
\overline f/1\in(\mathfrak nA_{\mathfrak n})^q,
\]
a contradiction. Hence \(f\notin\Phiop_q(I)\). This excludes every element of \(R\setminus I\) from some stage and independently proves
\[
\bigcap_{q\geq1}\Phiop_q(I)=I.
\]
\end{proof}

\begin{theorem}\label{thm:nontranslation-jet}
Let \(k\) be algebraically closed, let \(R=k[x_1,\ldots,x_n]\) with \(n\geq1\), and let \(q\geq1\). Then there is no ideal \(J\lhd R\) such that
\[
\Phiop_q(I)=I+J
\]
for every ideal \(I\lhd R\).
\end{theorem}

\begin{proof}
Corollary~\ref{cor:full-nullstellensatz} applied to \(I=(0)\) gives
\[
\sqrt{\Phiop_q(0)}=0.
\]
The polynomial ring \(R\) is reduced, so \(\Phiop_q(0)=0\). If a translating ideal \(J\) existed, taking \(I=(0)\) would give \(J=0\), and \(\Phiop_q\) would be the identity on \(\Id(R)\).

Let
\[
I=(x_1^{q+1},x_2,\ldots,x_n).
\]
Its affine zero set is the origin. If \(a\neq0\), at least one generator of \(I\) has nonzero value at \(a\), and the argument of equation~\eqref{eq:outside-factor} gives
\[
I+\m_a^q=R.
\]
At the origin, write \(\m_0=(x_1,\ldots,x_n)\). Every monomial of total degree \(q\) either equals \(x_1^q\) or contains one of \(x_2,\ldots,x_n\). Hence
\[
\m_0^q\subseteq(x_1^q,x_2,\ldots,x_n).
\]
The reverse generators satisfy \(x_1^q\in\m_0^q\) and \(x_i\in I\) for \(i\geq2\). Therefore
\[
I+\m_0^q=(x_1^q,x_2,\ldots,x_n).
\]
It follows that
\[
\Phiop_q(I)
=(x_1^q,x_2,\ldots,x_n)
\neq
(x_1^{q+1},x_2,\ldots,x_n)=I.
\]
This contradicts the identity conclusion. Hence no fixed translating ideal exists.
\end{proof}

\begin{theorem}\label{thm:univariate-formula}
Let \(k\) be algebraically closed, let \(R=k[t]\), let \(q\geq1\), and let
\[
g=c\prod_{i=1}^{s}(t-a_i)^{e_i},
\]
where \(c\in k^{\times}\), the elements \(a_i\in k\) are pairwise distinct, and the integers \(e_i\) are positive. Then
\[
\Phiop_q((g))
=
\left(\prod_{i=1}^{s}(t-a_i)^{\min\{e_i,q\}}\right).
\]
\end{theorem}

\begin{proof}
Let \(a\in k\setminus\{a_1,\ldots,a_s\}\). Then \(g(a)\neq0\), so equation~\eqref{eq:outside-factor} gives
\[
(g)+(t-a)^q=R.
\]
Such points do not affect the intersection.

Fix \(i\). In the principal ideal domain \(k[t]\), the sum of two principal ideals is generated by the greatest common divisor of their generators \cite[Chapter~1]{AtiyahMacdonald1969}. Since
\[
\prod_{j\neq i}(t-a_j)^{e_j}
\]
is relatively prime to \(t-a_i\),
\[
(g)+(t-a_i)^q
=
(t-a_i)^{\min\{e_i,q\}}.
\]
The resulting ideals are pairwise comaximal because the points \(a_i\) are distinct. Their intersection therefore equals their product by the Chinese remainder theorem \cite[Proposition~1.10]{AtiyahMacdonald1969}. Hence
\[
\begin{aligned}
\Phiop_q((g))
&=
\bigcap_{i=1}^{s}(t-a_i)^{\min\{e_i,q\}}\\
&=
\prod_{i=1}^{s}(t-a_i)^{\min\{e_i,q\}}\\
&=
\left(\prod_{i=1}^{s}(t-a_i)^{\min\{e_i,q\}}\right).
\end{aligned}
\]
\end{proof}

\begin{example}\label{ex:univariate}
Let \(k\) be algebraically closed and let
\[
I=(t^7(t-1)^3)\subseteq k[t].
\]
Then
\[
\Phiop_2(I)=(t^2(t-1)^2),
\qquad
\Phiop_4(I)=(t^4(t-1)^3),
\qquad
\Phiop_7(I)=I.
\]
These ideals have the same reduced zero set \(\{0,1\}\), while the observation order determines the retained multiplicity at each point.
\end{example}

\begin{theorem}\label{thm:boundedness-sharp}
Let \(k\) be an algebraically closed field, let \(R=k[x,y]\), and let \(c_1,c_2,\ldots\) be pairwise distinct elements of \(k\). Set
\[
X=\{(c_j,0):j\geq1\}
\]
and define
\[
\nu(c_j,0)=j.
\]
Then
\[
y\in\calI(X)
\qquad\text{and}\qquad
y\notin\sqrt{\Phiop_{X,\nu}(0)}.
\]
Consequently, the boundedness hypothesis in Theorem~\ref{thm:relative-nullstellensatz} cannot be omitted in general.
\end{theorem}

\begin{proof}
Every point of \(X\) has second coordinate zero, so
\[
y\in\calI(X).
\]
Assume that
\[
y\in\sqrt{\Phiop_{X,\nu}(0)}.
\]
Then there exists a positive integer \(N\) such that
\[
y^N\in\Phiop_{X,\nu}(0)
=
\bigcap_{j\geq1}(x-c_j,y)^j.
\]
Choose \(j>N\). The assumed membership gives
\[
y^N\in(x-c_j,y)^j.
\tag{6.3}\label{eq:sharp-membership}
\]
Consider the homomorphism
\[
\varepsilon_j:k[x,y]\longrightarrow k[y],
\qquad
\varepsilon_j(x)=c_j,
\qquad
\varepsilon_j(y)=y.
\]
The image of \((x-c_j,y)^j\) is contained in \((y^j)\). Applying \(\varepsilon_j\) to \eqref{eq:sharp-membership} gives
\[
y^N\in(y^j),
\]
which is impossible because \(N<j\). Thus no positive power of \(y\) belongs to \(\Phiop_{X,\nu}(0)\), proving the assertion.
\end{proof}

\begin{proposition}\label{prop:finite-computation}
Let \(k\) be a field, let \(R=k[x_1,\ldots,x_n]\), let \(s\geq1\), and let
\[
X=\{a_1,\ldots,a_s\}\subseteq k^n
\]
consist of distinct points, let \(\nu:X\to\bbN_{>0}\), and put \(q_i=\nu(a_i)\). Set
\[
J_{X,\nu}=\prod_{i=1}^{s}\m_{a_i}^{q_i}.
\]
Then, for every ideal \(I\lhd R\),
\[
\Phiop_{X,\nu}(I)
=
\bigcap_{i=1}^{s}\bigl(I+\m_{a_i}^{q_i}\bigr)
=
I+J_{X,\nu}.
\]
Moreover,
\[
J_{X,\nu}
=
\bigcap_{i=1}^{s}\m_{a_i}^{q_i}.
\]
Consequently, the restriction of a finite-observation nucleus to \(\Id(R)\) is translation by the fixed fat-point ideal \(J_{X,\nu}\).
\end{proposition}

\begin{proof}
For \(i\neq j\), the distinct maximal ideals \(\m_{a_i}\) and \(\m_{a_j}\) are comaximal. We verify that their positive powers are also comaximal. Choose
\[
u+v=1,
\qquad
u\in\m_{a_i},
\qquad
v\in\m_{a_j}.
\]
Expand
\[
1=(u+v)^{q_i+q_j-1}.
\]
In every term, either the exponent of \(u\) is at least \(q_i\), or the exponent of \(v\) is at least \(q_j\). Hence
\[
1\in\m_{a_i}^{q_i}+\m_{a_j}^{q_j}.
\]
Thus the ideals \(\m_{a_i}^{q_i}\) are pairwise comaximal. It follows that the ideals
\[
I+\m_{a_i}^{q_i}
\]
are pairwise comaximal as well.

For a finite pairwise comaximal family, the intersection equals the product by the Chinese remainder theorem \cite[Proposition~1.10]{AtiyahMacdonald1969}. Therefore
\[
\bigcap_{i=1}^{s}(I+\m_{a_i}^{q_i})
=
\prod_{i=1}^{s}(I+\m_{a_i}^{q_i}).
\tag{6.4}\label{eq:comax-product}
\]
Expanding the product on the right, every term that contains a factor from \(I\) belongs to \(I\), while the unique term containing no such factor belongs to
\[
\prod_{i=1}^{s}\m_{a_i}^{q_i}=J_{X,\nu}.
\]
Hence
\[
\prod_{i=1}^{s}(I+\m_{a_i}^{q_i})
\subseteq
I+J_{X,\nu}.
\]
Conversely, both \(I\) and \(J_{X,\nu}\) are contained in every factor \(I+\m_{a_i}^{q_i}\). Thus
\[
I+J_{X,\nu}
\subseteq
\bigcap_{i=1}^{s}(I+\m_{a_i}^{q_i}).
\]
Together with \eqref{eq:comax-product}, this proves
\[
\Phiop_{X,\nu}(I)=I+J_{X,\nu}.
\]
Applying the same pairwise-comaximal intersection-product identity to the ideals \(\m_{a_i}^{q_i}\) gives
\[
J_{X,\nu}=\bigcap_{i=1}^{s}\m_{a_i}^{q_i}.
\]
\end{proof}

\medskip

\subsection*{Acknowledgements}
The author gratefully acknowledges the referee for a careful and constructive reading of the manuscript and for the detailed comments and suggestions. These recommendations led to substantial improvements in the organization, mathematical formulation, bibliographic positioning, and exposition of the paper. In particular, the referee's suggestions concerning the foundational definitions, the proximity literature, the terminology for prime ideals, and the presentation of several statements have helped clarify the mathematical scope and improve the overall coherence of the paper.

\subsection*{Funding}

None.

\subsection*{Declaration of Competing Interest}

The author declares no competing interest.

\subsection*{Data Availability}

No data were used in this research.

\subsection*{Declaration on AI usage}
All ideas, as well as credit or blame for the mathematical arguments, belong strictly to the human author. AI was used solely for proofreading the manuscript and verifying the correctness of the calculations.

\end{document}